\documentclass[12pt,a4paper]{article}
\usepackage[english]{babel}
\usepackage{amssymb}
\usepackage[T1]{fontenc}
\usepackage{amsmath}
\newtheorem{Th}{Theorem}

\newtheorem{Lem}{Lemma}

\newcommand{\A}{\mathcal{A}}
\newcommand{\F}{\mathcal{F}}
\newcommand{\R}{\mathbb{R}}
\newcommand{\Z}{\mathbb{Z}}
\renewcommand{\H}{\mathbb{H}}
\newcommand{\N}{\mathbb{N}}

\newcommand{\B}{\mathcal{B}}

\newcommand{\C}{\mathcal{C}}

\newcommand{\Q}{\mathcal{Q}}

\renewcommand{\P}{\mathbb{P}}
\newcommand{\ds}{\displaystyle}
\newcounter{tictac}

\def\1{\,\rlap{\mbox{\small\rm 1}}\kern.15em 1}
\def\ind#1{\1_{#1}}
\def\build#1_#2^#3{\mathrel{\mathop{\kern 0pt#1}\limits_{#2}^{#3}}}
\def\tend#1#2{\build\hbox to 12mm{\rightarrowfill}_{#1\rightarrow #2}^{ }}

\def\converge#1#2#3#4{\build\hbox to
#1mm{\rightarrowfill}_{#2\rightarrow #3}^{\hbox{\scriptsize #4}}}
\begin{document}
\title{Invariance principles for standard-normalized and self-normalized random fields}
\author{$\textrm{M. El Machkouri}^\ast$, L. Ouchti\thanks{LMRS, UMR CNRS 6085, Universit\'e de Rouen, Site Colbert, 76821 Mont Saint Aignan, France. 
E-mail : mohamed.elmachkouri@univ-rouen.fr, lahcen.ouchti@univ-rouen.fr}}
\maketitle
\begin{abstract}
We investigate the invariance principle for set-indexed partial
sums of a stationary field $(X_{k})_{k\in\Z^{d}}$ of
martingale-difference or independent random variables under
standard-normalization or self-normalization respectively.\\
\\
{\em AMS Classifications (2000)} : 60 F 05, 60 F 17, 60 G 60
{\em Short title} : Invariance principles for random fields
\end{abstract}

\section{Introduction}
Let $(X_{k})_{k\in\Z^{d}}$ be a stationary field of real-valued
random variables defined on a given probability space $(\Omega,
\F, \P)$. If $\A$ is a collection of Borel subsets of $[0,1]^{d}$,
define the smoothed partial sum process $\{S_{n}(A)\,;\,A\in\A\}$
by
\begin{equation}\label{process}
S_{n}(A)=\ds{\sum_{i\in\{1,...,n\}^{d}}}\,\lambda(nA\cap
R_{i})X_{i}
\end{equation}
where $R_{i}=]i_{1}-1,i_{1}]\times...\times]i_{d}-1,i_{d}]$ is the
unit cube with upper corner at $i$ and $\lambda$ is the Lebesgue
measure on $\R^{d}$. We equip the collection $\A$ with the
pseudo-metric $\rho$ defined for any $A,B$ in $\A$ by
$\rho(A,B)=\sqrt{\lambda(A\Delta B)}$. To measure the size of $\A$
one considers the metric entropy: denote by
$H(\A,\rho,\varepsilon)$ the logarithm of the smallest number
$N(\A,\rho,\varepsilon)$ of open balls of radius $\varepsilon$
with respect to $\rho$ which form a covering of $\A$. The function
$H(\A, \rho, .)$ is the entropy of the class $\A$. A more strict
tool is the metric entropy with inclusion: assume that $\A$ is
totally bounded with inclusion i.e. for each positive
$\varepsilon$ there exists a finite collection $\A(\varepsilon)$
of Borel subsets of $[0,1]^{d}$ such that for any $A\in\A$, there
exist $A^{-}$ and $A^{+}$ in $\A(\varepsilon)$ with
$A^{-}\subseteq A\subseteq A^{+}$ and
$\rho(A^{-},A^{+})\leq\varepsilon$. Denote by
$\mathbb{H}(\A,\rho,\varepsilon)$ the logarithm of the cardinality
of the smallest collection $\A(\varepsilon)$. The function
$\mathbb{H}(\A,\rho,.)$ is the entropy with inclusion (or
bracketing entropy) of the class $\A$. Let $C(\A)$ be the space of
continuous real functions on $\A$, equipped with the norm $\|
.\|_{\A}$ defined by
$$
\| f\|_{\A}=\sup_{A\in\A}\vert f(A)\vert.
$$
A standard Brownian motion indexed by $\A$ is a mean zero Gaussian
process $W$ with sample paths in $C(\A)$ and
Cov(W(A),W(B))$=\lambda(A\cap B)$. From Dudley \cite{Dudley} we
know that such a process exists if
\begin{equation}\label{entrop-metriq1}
\int_{0}^{1}\sqrt{H(\A,\rho,\varepsilon)}\,d\varepsilon<+\infty.
\end{equation}
Since $H(\A,\rho,.)\leq\mathbb{H}(\A,\rho,.)$, the standard
Brownian motion $W$ is well defined if
\begin{equation}\label{entrop-metriq2}
\int_{0}^{1}\sqrt{\mathbb{H}(\A,\rho,\varepsilon)}\,d\varepsilon<+\infty.
\end{equation}
For any probability measure $m$ defined on $[0,1]^{d}$ equipped
with its Borel $\sigma$-algebra, we define the pseudo-metric
$\rho_{m}$ by $\rho_{m}=\sqrt{m(A\Delta B)}$ for any $A$ and $B$
in $\A$. For any positive $\varepsilon>0$, we denote
$N(\A,\varepsilon)=\sup_{m}N(\A,\rho_{m},\varepsilon)$ and we say
that the collection $\A$ has uniformly integrable entropy if
\begin{equation}\label{uniform-entropy}
\int_{0}^{1}\sqrt{\log N(\A,\varepsilon)}d\varepsilon<+\infty.
\end{equation}
We say that the (classical) invariance principle or functional
central limit theorem (FCLT) holds if the sequence
$\{n^{-d/2}S_{n}(A)\,;\,A\in\A\}$ converges in distribution to an
$\A$-indexed Brownian motion in the space $C(\A)$. The first weak
convergence results for $\Q_{d}$-indexed partial sum processes
were established for i.i.d. random fields and for the collection
$\Q_{d}$ of lower-left quadrants in $[0,1]^{d}$, that is to say
the collection
$\{[0,t_{1}]\times\ldots\times[0,t_{d}]\,;\,(t_{1},\ldots,t_{d})\in[0,1]^{d}\}$.
They were proved by Wichura \cite{Wichu} under a finite variance
condition and earlier by Kuelbs \cite{Kuelbs} under additional
moment restrictions. When the dimension $d$ is reduced to one,
these results coincide with the original invariance principle of
Donsker \cite{Donsker}. In 1983, Pyke \cite{Pyke} derived a weak
convergence result for the process $\{S_{n}(A)\,;\,A\in\A\}$ for
i.i.d. random fields provided that the collection $\A$ satisfies
the bracketing entropy condition ($\ref{entrop-metriq2}$).
However, his result required moment conditions which depend on the
size of the collection $\A$. Bass \cite{Bass} and simultaneously
Alexander and Pyke \cite{Alex-Pyke} extended Pyke's result to
i.i.d. random fields
with finite variance. More precisely, the following result is proved.\\
\vspace{-0.25cm}
\\
\textbf{Theorem A\,(Bass (1985), Alexander and Pyke (1986))} {\em
Let $(X_{k})_{k\in\Z^{d}}$ be a stationary field of independent
real random variables with zero mean and finite variance. If $\A$
is a collection of regular Borel subsets of $[0,1]^{d}$ which
satisfies Assumption $(\ref{entrop-metriq2})$ then the sequence of
processes $\{n^{-d/2}S_{n}(A);A\in\A\}$ converge in distribution
to $\sqrt{E(X_{0}^{2})}W$ where $W$ is a
standard Brownian motion indexed by $\A$.}\\
\vspace{-0.25cm}
\\
Unfortunately, the bracketing condition ($\ref{entrop-metriq2}$)
is not automatically fulfilled in the important case of $\A$ being
a Vapnik-Chervonenkis class of sets. Ziegler \cite{Ziegler} has
covered this case by proving (among other results) that the FCLT
of Bass, Alexander and Pyke (i.e. Theorem A) still holds for
classes of sets which satisfy the uniformly integrable entropy
condition ($\ref{uniform-entropy}$). Recently, Dedecker
\cite{JD-tflc} gave an $L^{\infty}$-projective criterion for the
process $\{n^{-d/2}S_{n}(A)\,;\,A\in\A\}$ to converge to a mixture
of $\A$-indexed Brownian motions when the collection $\A$
satisfies only the entropy condition ($\ref{entrop-metriq1}$) of
Dudley. This new criterion is valid for martingale-difference
bounded random fields and provides a new criterion for non-uniform
$\phi$-mixing bounded random fields. In the unbounded case, using
the chaining method of Bass \cite{Bass} and establishing Bernstein
type inequalities, Dedecker proved also the FCLT for the partial
sum $\{S_{n}(A)\,;\,A\in\A\}$ of non-uniform $\phi$-mixing random
fields provided that the collection $\A$ satisfies the more strict
entropy condition with inclusion ($\ref{entrop-metriq2}$) and
under both finite fourth moments and a polynomial decay of the
mixing coefficients. In a previous work (see \cite{EM-Volny}), it
is shown that the FCLT may be not valid for $p$-integrable ($0\leq
p<+\infty$) martingale-difference random fields. More precisely, the following result is established.\\
\vspace{-0.25cm} \\
\textbf{Theorem B\,(El Machkouri, Voln\'y, 2002)} {\em Let
$(\Omega,\F,\mu,T)$ be an ergodic dynamical system with positive
entropy where $\Omega$ is a Lebesgue space, $\mu$ is a probability
measure and $T$ is a $\Z^{d}$-action. For any nonnegative real
$p$, there exist a real function $f\in L^{p}(\Omega)$ and a
collection $\A$ of regular Borel subsets of $[0,1]^{d}$ such that
\begin{itemize}
\item For any $k$ in $\Z^{d}$, $E\left(f\circ T^{k}\vert\sigma(f\circ
    T^{i}\,;\,i\neq k)\right)=0$. We say that the random field
    $(f\circ T^{k})_{k\in\Z^{d}}$ is a strong martingale-difference
random field.
\item The collection $\A$ satisfies the entropy condition with
inclusion $(\ref{entrop-metriq2})$.
\item The partial sum process $\{n^{-d/2}S_{n}(f,A)\,;\,A\in\A\}$ is not tight in the space
$C(\A)$
\end{itemize}
where}
$$
S_{n}(f,A):=\sum_{i\in\{1,...,n\}^{d}}\lambda(nA\cap R_{i})f\circ
T^{i}.
$$
The above theorem shows that not only Dedecker's FCLT for bounded
random fields (see \cite{JD-tflc}) cannot be extended to
$p$-integrable ($0\leq p<+\infty$) random fields but also it lays
emphasis on that Bass, Alexander and Pyke's result for i.i.d.
random fields (Theorem A) cannot hold for martingale-difference
random fields without additional assumptions. Recently, El
Machkouri \cite{MEM-KK} has shown that the FCLT still holds for
unbounded random fields which satisfy both a finite exponential
moment condition and a projective criterion similar to Dedecker's
one. All these results put on light that the moment assumption on
the random field is very primordial in the FCLT question for
random fields indexed by large classes of sets.\\
In the present work, we give a positive answer to the validity of
the FCLT for square-integrable martingale-difference random fields
which conditional variances are bounded almost surely (cf. Theorem
$\ref{fclt-am}$). Next, we consider self-normalized i.i.d. random
fields, more precisely, we investigate the validity of the FCLT
when the stationary random field $(X_{k})_{k\in\Z^{d}}$ is assumed
to be independent and the classical normalization $n^{d/2}$ is
replaced by $U_{n}$ defined by ($\ref{Un}$) (cf. Theorem
$\ref{fclt}$). From a statistical point of view, the
self-normalization is natural and several articles in the
literature are devoted to limit theorems for self-normalized
sequences $(X_{k})_{k\in\Z}$ of independent random variables with
statistical applications. Logan et al.
\cite{Logan-Mallows-Rice-Shepp} investigate the various possible
limit distributions of self-normalized sums. Gin\'e et al.
\cite{Gine-Gotze-Mason} prove that
$\sum_{i=1}^{n}X_{i}/\sqrt{\sum_{i=1}^{n}X_{i}^{2}}$ converges to
the Gaussian standard distribution if and only if $X_{1}$ is in
the domain of attraction of the normal distribution (the symmetric
case was previously treated by Griffin and Mason
\cite{Griffin-Mason}). Egorov \cite{Egorov} investigates the non
identically distributed case. Large deviations are investigated in
Shao \cite{Shao97} without moment conditions. Ra\v ckausksas and
Suquet \cite{Rackauskas-Suquet} gives invariance principles for
various partial sums processes under self-normalization in
$\C([0,1])$ and in the stronger topological framework of H\"older
spaces. Our Theorem $\ref{fclt}$ below improves on Ra\v ckauskas
and Suquet's result in $\C([0,1])$.
\section{Main results}
By a stationary real random field we mean any family
$(X_{k})_{k\in\Z^{d}}$ of real-valued random variables defined on
a probability space $(\Omega, \F, \P)$ such that for any
$(k,n)\in\Z^{d}\times\N^{\ast}$ and any
$(i_{1},...,i_{n})\in(\Z^{d})^n$, the random vectors
$(X_{i_{1}},...,X_{i_{n}})$ and $(X_{i_{1}+k},...,X_{i_{n}+k})$
have the same law.\\
On the lattice $\Z^{d}$ we define the lexicographic order as
follows: if $i=(i_{1},...,i_{d})$ and $j=(j_{1},...,j_{d})$ are
distinct elements of $\Z^{d}$, the notation $i<_{lex}j$ means that
either $i_{1}<j_{1}$ or for some $p$ in $\{2,3,...,d\}$,
$i_{p}<j_{p}$ and $i_{q}=j_{q}$ for $1\leq q<p$. A real random
field $(X_{k})_{k\in\Z^{d}}$ is said to be a martingale-difference
random field if it satisfies the following condition: for any $m$
in $\Z^{d}$, $E\left(X_{m}\vert\F_{m}\right)=0$ a.s. where
$\F_{m}$ is the $\sigma$-algebra generated by the random variables
$X_{k},\,k<_{lex}m$. Our first result is the following.
\begin{Th}\label{fclt-am}
Let $(X_{k})_{k\in\Z^{d}}$ be a stationary field of
martingale-difference random variables with finite variance such
that $E(X_{0}^2\vert\F_{0})$ is bounded almost surely and let $\A$
be a collection of regular Borel subsets of $[0,1]^{d}$ satisfying
the condition $(\ref{entrop-metriq2})$. Then the sequence
$\{n^{-d/2}S_{n}(A);\, A\in\A\}$ converges weakly in $C(\A)$ to
$\sqrt{E(X_{0}^2)}W$ where $W$ is the standard Brownian motion
indexed by $\A$.
\end{Th}
Comparing Theorem $\ref{fclt-am}$ and Theorem B in section 1, one
can notice that the conditional variance
$E\left(X_{0}^2\vert\F_{0}\right)$ is primordial in the invariance
principle problem for martingale-difference random fields. More
generally, the conditional variance for martingales is known to
play an important role in modern
martingale limit theory (see Hall and Heyde \cite{Hall-Heyde}).\\
For any integer $n\geq 1$, we define
\begin{equation}\label{Un}
U_{n}^{2}=\sum_{i\in\Lambda_{n}}X_{i}^{2}
\end{equation}
where $\Lambda_{n}=\{1,...,n\}^{d}$. We say that $X_{0}$ belongs
to the domain of attraction of the normal distribution (and we
denote $X_{0}\in DAN$) if there exists a norming sequence $b_{n}$
of real numbers such that $b_{n}^{-1}S_{\Lambda_{n}}$ converges in
distribution to a standard normal law. We should recall that if
$X_{0}\in DAN$ then $\|X_{0}\|_{p}<\infty$ for any $0<p<2$ and
that constants $b_{n}$ have the form $b_{n}=n^{d/2}l(n)$ for some
function $l$ slowly varying at infinity. Moreover, for each
$\tau>0$, we have
\begin{equation}\label{standard}
\lim_{n\to\infty}n^d EX_{0,n}=0,\,\,\lim_{n\to\infty}n^{d}\P(\vert
X_{0}\vert\geq\tau b_{n})=0\,\,\textrm{and}\,
\lim_{n\to\infty}b_{n}^{-2}n^d E(X_{0,n}^2)=1
\end{equation}
where $X_{0,n}=X_{0}\ind{\vert X_{0}\vert<\tau b_{n}}$ (see for
instance Araujo and Gin\'e \cite{Araujo-Gine}). Note also that
$X_{0}\in DAN$ implies (Raikov's theorem) that
\begin{equation}\label{Raikov}
\frac{1}{b_{n}^2}\sum_{i\in\Lambda_{n}}X_{i}^2\converge{12}{n}{\infty}{$\P$}1.
\end{equation}
\begin{Th}\label{fclt}
Let $(X_{k})_{k\in\Z^{d}}$ be a field of i.i.d. centered random
variables and let $\A$ be a collection of regular Borel subsets of
$[0,1]^{d}$ satisfying the condition $(\ref{entrop-metriq2})$.
Then $X_{0}\in DAN$ if and only if the sequence
$\{U_{n}^{-1}S_{n}(A);\, A\in\A\}$ converges weakly in $C(\A)$ to
the standard Brownian motion $W$.
\end{Th}
Let us remark that the necessity of $X_{0}\in DAN$ in Theorem
$\ref{fclt}$ follows from Gin\'e et al. (\cite{Gine-Gotze-Mason},
Theorem 3.3). Our result contrasts with the invariance principle
established by Bass and Alexander and Pyke (cf. Theorem A in
section 1) where square integrable random variables are required.
We do not know if Theorem $\ref{fclt}$ still hold if one replace
the condition ($\ref{entrop-metriq2}$) by condition
($\ref{entrop-metriq1}$). However, our next result is a
counter-example which shows that Theorem A in section 1 does not
hold when the condition ($\ref{entrop-metriq2}$) is replaced by
condition ($\ref{entrop-metriq1}$).
\begin{Th}\label{counter}
For any positive real number $p$, there exist a stationary field
$(X_{k})_{k\in\Z^{d}}$ of independent, symmetric and
$p$-integrable real random variables and a collection $\A$ of
regular Borel subsets of $[0,1]^{d}$ which satisfies the condition
$(\ref{entrop-metriq1})$ such that the partial sum process
$\{n^{-d/2}S_{n}(A)\,;\,A\in\A\}$ do not be tight in the space
$C(\A)$.
\end{Th}
Note that Dudley and Strassen \cite{Dudley-Strassen} have built a
sequence of i.i.d. random variables $X_{n}$ with values in the
space of continuous functions on $[0,1]$ such that $E(X_{1}(t))=0$
and the finite dimensional marginals of
$Z_{n}(t)=n^{-1/2}\sum_{i=1}^{n}X_{i}(t)$ converge to that of a
Gaussian process $Z$. It was shown that this process $Z$ has a
version with almost sure continuous sample paths and that the
process $Z_{n}(t)$ is not tight for the topology of the uniform
metric. However, contrary to our example, one can check that the
limiting process $Z$ does not satisfy the Dudley's entropy
condition ($\ref{entrop-metriq1}$) for the intrinsic distance
$\rho(s,t)=\|Z(s)-Z(t)\|_{2}$. In fact, it is well known that the
condition ($\ref{entrop-metriq1}$) is sufficient for Gaussian
processes to have a version with almost sure continuous sample
paths but it falls to be necessary (see van der Vaart and Wellner
\cite{van-der-Vaart-Wellner}, p. 445).
\section{Proofs}
Recall that a Young function $\psi$ is a real convex nondecreasing
function defined on $\R^{+}$ which satisfies $\psi(0)=0$. We
define the Orlicz space $L_{\psi}$ as the space of real random
variables $Z$ defined on the probability space $(\Omega, \F, \P)$
such that $E[\psi(\vert Z\vert/c)]<+\infty$ for some $c>0$. The
Orlicz space $L_{\psi}$ equipped with the so-called Luxemburg norm
$\| . \|_{\psi}$ defined for any real random variable $Z$ by
$$
\| Z\|_{\psi}=\inf\{\,c>0\,;\,E[\psi(\vert Z\vert/c)]\leq 1\,\}
$$
is a Banach space. For more about Young functions and Orlicz
spaces one can refer to Krasnosel'skii and Rutickii \cite{K-R}.
Let $\psi_{1},\psi_{2}:\R^{+}\to\R$ be the Young functions defined
by $\psi_{1}(x)=\exp(x)-1$ and $\psi_{2}(x)=\exp(x^2)-1$ for any
$x\in\R^{+}$. We need the following lemma which is of independent
interest.
\begin{Lem}\label{Lem1}
Let $(\theta_{i})_{i\in\Z^d}$ be an arbitrary field of random
variables and let $\mathcal{H}_{i}$ denote the $\sigma$-algebra
generated by the random variables
$\theta_{j},\,j<_{lex}i,\,i\in\Z^d$. Let also
$0\leq\alpha\leq\beta\leq 1$ and $0<\tau\leq 1$ be fixed and let
$(c_{n})_{n\geq 1}$ be a sequence of real numbers. For any integer
$n\geq 1$ and any Borel subset $A$ of $[0,1]^d$, denote
$$
\theta_{i}(n,\alpha,\beta)=\theta_{i}\ind{{\alpha\tau
c_{n}\leq\vert\theta_{i}\vert<\beta\tau c_{n}}}
$$
and
$$
\Theta_{n}(A,\alpha,\beta)=\frac{1}{c_{n}}\sum_{i\in\Lambda_{n}}\lambda(nA\cap
R_{i})[\theta_{i}(n,\alpha,\beta)-E\left(\theta_{i}(n,\alpha,\beta)\vert\mathcal{H}_{i})\right].
$$
Assume also that there exists $C>0$ such that for any integer
$n\geq 1$ and any $i$ in $\Z^d$,
\begin{equation}\label{hypothese}
\frac{n^d}{c_{n}^2}E\left(\theta_{i}^2\ind{\vert\theta_{i}\vert<
c_{n}}\vert\mathcal{H}_{i}\right)\leq C.
\end{equation}
If $\mathcal{G}_{1},\mathcal{G}_{2}$ are finite collections of
Borel subsets of $[0,1]^d$ then
$$
\bigg\|\max_{(A,B)\in\mathcal{G}}\big\vert\Theta_{n}(A,\alpha,\beta)-\Theta_{n}(B,\alpha,\beta)\big\vert\bigg\|_{\psi_{1}}
\hspace{-0.2cm}\leq
K[\beta\,\tau\,\psi_{1}^{-1}(\vert\mathcal{G}\vert)+
\max_{(A,B)\in\mathcal{G}}\rho(A,B)\,\psi_{2}^{-1}(\vert\mathcal{G}\vert)]
$$
where $\mathcal{G}=\mathcal{G}_{1}\times\mathcal{G}_{2}$,
$\vert\mathcal{G}\vert$ is the cardinal of $\mathcal{G}$ and $K>0$
is a universal constant.
\end{Lem}
{\em Proof of Lemma $\ref{Lem1}$}. Consider the field of
martingale-difference random variables
$Y_{i}(n,\alpha,\beta),\,i\in\Lambda_{n}$ defined by
$$
Y_{i}(n,\alpha,\beta)= \frac{1}{c_{n}}(\lambda(nA\cap
R_{i})-\lambda(nB\cap
R_{i}))[\theta_{i}(n,\alpha,\beta)-E\left(\theta_{i}(n,\alpha,\beta)\vert\mathcal{H}_{i}\right)]
$$
and note that $\vert Y_{i}(n,\alpha,\beta)\vert\leq 2\beta\tau$.
Using ($\ref{hypothese}$) and keeping in mind that $\tau$ and
$\beta$ are less than 1, there exists a universal constant $C>0$
such that
$$
\sum_{i\in\Lambda_{n}}E\left(Y_{i}(n,\alpha,\beta)^2\vert\mathcal{H}_{i}\right)\leq
4C\max_{(A,B)\in\mathcal{G}}\rho^{2}(A,B).
$$
Noting that
$\Theta_{n}(A,\alpha,\beta)-\Theta_{n}(B,\alpha,\beta)=\sum_{i\in\Lambda_{n}}
Y_{i}(n,\alpha,\beta)$ and applying Theorem 1.2A in de la Pena
\cite{Pena}, we derive the following Bernstein inequality
$$
\P\left(\big\vert
\Theta_{n}(A,\alpha,\beta)-\Theta_{n}(B,\alpha,\beta)\big\vert>x\right)
\leq
2\exp\left(\frac{-x^2}{8C\max_{(A,B)\in\mathcal{G}}\rho^2(A,B)+4\beta\tau
x}\right).
$$
The proof is completed by using Lemma 2.2.10 in van der Vaart and
Wellner \cite{van-der-Vaart-Wellner}.
\subsection{Proof of Theorem $\bold{\ref{fclt-am}}$}
$\small{\textbf{\quad a) Tightness}}$\\
\vspace{-0.2cm}
\\
It suffices to prove that for any $x>0$
\begin{equation}\label{tension-am}
\lim_{\delta\to0}\limsup_{n\to+\infty}
\P\left(\sup_{\substack{A,B\in\A \\
\rho(A,B)<\delta}}\big\vert
n^{-d/2}S_{n}(A)-n^{-d/2}S_{n}(B)\big\vert>x\right)=0.
\end{equation}
In the sequel, we write $\H(x)$ for $\H(\A,\rho,x)$. Let
$\delta>0$ be fixed, denote $\tau=\delta/\sqrt{\H(\delta/2)}>0$
and assume (without loss of generality) that $\tau\leq 1$. Let
$i\in\Z^d$, since $X_{i}$ is a martingale-difference random
variable, we have
$X_{i}=X_{i,n}-E(X_{i,n}\vert\F_{i})+\overline{X}_{i,n}-E(\overline{X}_{i,n}\vert\F_{i})$
where $X_{i,n}=X_{i}\ind{{\vert X_{i}\vert<\tau n^{d/2}}}$ and
$\overline{X}_{i,n}=X_{i}-X_{i,n}$, hence it follows
$$
\P\left(\sup_{\substack{A,B\in\A \\
\rho(A,B)<\delta}}\big\vert
n^{-d/2}S_{n}(A)-n^{-d/2}S_{n}(B)\big\vert>x\right)\leq
E_{1}+E_{2}
$$
where
\begin{align*}
E_{1}&=\P\left(\sup_{{\substack{A,B\in\A \\
\rho(A,B)<\delta}}}\bigg\vert\sum_{i\in\Lambda_{n}}(\lambda(nA\cap R_{i})-\lambda(nB\cap R_{i}))[X_{i,n}-E\left(X_{i,n}\vert\F_{i}\right)]\bigg\vert>xn^{d/2}/2\right)\\
E_{2}&=n^d\P\left(\vert X_{0}\vert\geq\tau
n^{d/2}\right)\tend{n}{+\infty}0\qquad\textrm{(since $X_{0}\in
L^2$)}.
\end{align*}
We are going to control $E_{1}$. Now, for any constants
$0\leq\alpha\leq\beta\leq 1$ define
$X_{i}(n,\alpha,\beta)=X_{i}\ind{\alpha\tau n^{d/2}\leq\vert
X_{i}\vert<\beta\tau n^{d/2}}$ and
$$
Z_{n}(A,\alpha,\beta)=\frac{1}{n^{d/2}}\sum_{i\in\Lambda_{n}}\lambda(nA\cap
R_{i})[X_{i}(n,\alpha,\beta)-E\left(X_{i}(n,\alpha,\beta)\vert\F_{i}\right)].
$$
One can notice that
$$
E_{1}\leq\frac{2}{x}E\left(\sup_{\substack{A,B\in\A \\
\rho(A,B)<\delta}}\big\vert
Z_{n}(A,0,1)-Z_{n}(B,0,1)\big\vert\right).
$$
Let $\delta_{k}=2^{-k}\delta$. If $A$ and $B$ are any sets in
$\A$, there exists sets $A_{k},A_{k}^{+},B_{k},B_{k}^{+}$ in the
finite class $\A(\delta_{k})$ such that $A_{k}\subset A\subset
A_{k}^{+}$ and $\rho(A_{k},A_{k}^{+})\leq\delta_{k}$, and
similarly for $B,B_{k},B_{k}^{+}$. Let $(a_{k})_{k\in\N}$ be a
sequence of positive numbers decreasing to zero sucht that
$a_{0}=1$. Following the chaining method initiated by Bass
\cite{Bass}, we write
\begin{align*}
Z_{n}(A,0,1)-Z_{n}(A_{0},0,1)&=\sum_{k=0}^{+\infty}Z_{n}(A_{k+1},0,a_{k})-Z_{n}(A_{k},0,a_{k})\\
&\quad
+\sum_{k=1}^{+\infty}Z_{n}(A,a_{k},a_{k-1})-Z_{n}(A_{k},a_{k},a_{k-1}).
\end{align*}
So, we have $ \frac{x}{2}E_{1}\leq F_{1}+F_{2}+F_{3}$ where
\begin{align*}
F_{1}&=E\left(\max_{\substack{A_{0},\,B_{0}\in\A(\delta_{0}) \\
\rho(A_{0},B_{0})\leq 3\delta_{0}}}\big\vert Z_{n}(A_{0},0,1)-Z_{n}(B_{0},0,1)\big\vert\right)\\
F_{2}&=2\sum_{k=0}^{+\infty}
E\left(\max_{\substack{A_{k}\in\A(\delta_{k}),\,A_{k+1}\in\A(\delta_{k+1})\\
\rho(A_{k},A_{k+1})\leq 2\delta_{k}}}
\big\vert Z_{n}(A_{k+1},0,a_{k})-Z_{n}(A_{k},0,a_{k})\big\vert\right)\\
F_{3}&=2\sum_{k=1}^{+\infty}
E\left(\max_{\substack{A_{k},\,A_{k}^{+}\in\A(\delta_{k})\\
\rho(A_{k},A_{k}^{+})\leq\delta_{k}}}\,\sup_{A_{k}\subset A\subset
A_{k}^{+}}\big\vert
Z_{n}(A,a_{k},a_{k-1})-Z_{n}(A_{k},a_{k},a_{k-1})\big\vert\right)
\end{align*}
In the sequel, we denote by $K$ any universal positive constant.
Applying Lemma $\ref{Lem1}$ with $c_{n}=n^{d/2}$, we derive
\begin{equation}\label{F1}
F_{1}\leq
K\,\left(\tau\H(\delta_{0})+\delta_{0}\sqrt{\H(\delta_{0})}\right),
\end{equation}
similarly
\begin{equation}\label{F2}
F_{2}\leq
K\sum_{k=0}^{+\infty}(a_{k}\tau\H(\delta_{k+1})+\delta_{k}\sqrt{\H(\delta_{k+1})}).
\end{equation}
Now, we are going to control the last term $F_{3}$. For any Borel
subset $A$ of $[0,1]^d$, we denote
$$
\widetilde{Z}_{n}(A,a_{k},a_{k-1})=\frac{1}{n^{d/2}}\sum_{i\in\Lambda_{n}}\lambda(nA\cap
R_{i})[\vert X_{i}(n,a_{k},a_{k-1})\vert-E\left(\vert
X_{i}(n,a_{k},a_{k-1})\vert\vert\F_{i}\right)].
$$
One can check that
\begin{align*}
&\sup_{A_{k}\subset A\subset A_{k}^{+}}\vert Z_{n}(A,a_{k},a_{k-1})-Z_{n}(A_{k},a_{k},a_{k-1})\vert\\
&\qquad\leq
\frac{1}{n^{d/2}}\sum_{i\in\Lambda_{n}}(\lambda(nA_{k}^{+}\cap
R_{i})-\lambda(nA_{k}\cap R_{i}))[\vert
X_{i}(n,a_{k},a_{k-1})\vert-E\left(\vert
X_{i}(n,a_{k},a_{k-1})\vert\vert\F_{i}\right)]\\
&\qquad\qquad+\frac{2}{n^{d/2}}\sum_{i\in\Lambda_{n}}(\lambda(nA_{k}^{+}\cap
R_{i})-\lambda(nA_{k}\cap R_{i}))E\left(\vert
X_{i}(n,a_{k},a_{k-1})\vert\vert\F_{i}\right)\\
&\qquad=\widetilde{Z}_{n}(A_{k}^{+},a_{k},a_{k-1})-\widetilde{Z}_{n}(A_{k},a_{k},a_{k-1})\\
&\qquad\qquad+\frac{2}{n^{d/2}}\sum_{i\in\Lambda_{n}}\lambda(n\left(A_{k}^{+}\backslash
A_{k}\right)\cap R_{i})E\left(\vert
X_{i}(n,a_{k},a_{k-1})\vert\vert\F_{i}\right)
\end{align*}
Recall that by assumption we have $E(X_{i}^2\vert\F_{i})\leq C$
for some $C>0$. So, using Lemma $\ref{Lem1}$, it follows
$$
\bigg\|\max_{A_{k},A_{k}^{+}\in\A(\delta_{k})}\big\vert\widetilde{Z}_{n}(A_{k}^{+},a_{k},a_{k-1})-\widetilde{Z}_{n}(A_{k},a_{k},a_{k-1})\big\vert\bigg\|_{\psi_{1}}\hspace{-0.2cm}\leq
K(a_{k-1}\tau\H(\delta_{k})+\delta_{k}\sqrt{\H(\delta_{k})}).
$$
Moreover, one can check that
$$
E\left(\vert X_{i}(n,a_{k},a_{k-1})\vert\vert\F_{i}\right)
\leq\frac{E\left(X_{i}^2\vert\F_{i}\right)}{a_{k}\tau
n^{d/2}}\leq\frac{C}{a_{k}\tau n^{d/2}}.
$$
Consequently, we obtain
\begin{equation}\label{F3}
F_{3}\leq
K\left(\sum_{k=1}^{+\infty}a_{k-1}\tau\H(\delta_{k})+\delta_{k}\sqrt{\H(\delta_{k})}+\frac{\delta_{k}^2}{\tau
a_{k}}\right)
\end{equation}
Now, we choose $a_{k}=\delta_{k}/(\tau\sqrt{\H(\delta_{k+1})})$
for all $k\in\N$ (note that $a_{0}=1$), hence, we obtain the
following estimations:
\begin{align*}
F_{1}&\leq K\,\delta\sqrt{\H(\delta/2)}\\
F_{2}&\leq K\sum_{k=0}^{+\infty}\delta_{k}\sqrt{\H(\delta_{k+1})}\\
F_{3}&\leq
K\sum_{k=1}^{+\infty}\delta_{k-1}\sqrt{\H(\delta_{k+1})}
\end{align*}
Now, recall that $\frac{2}{x}E_{1}\leq F_{1}+F_{2}+F_{3}$ and keep
in mind that the entropy condition ($\ref{entrop-metriq2}$) holds
then
$$
\limsup_{n\to\infty}\frac{2}{x}E_{1}\leq
K\sum_{k=1}^{+\infty}\delta_{k+1}\sqrt{\H(\delta_{k})}\leq
K\int_{0}^{\delta}\sqrt{\H(x)}dx\tend{\delta}{0}0.
$$
Finally, the condition ($\ref{tension-am}$) holds and the sequence
$\{n^{-d/2}S_{n}(A)\,;\,A\in\A\}$ is tight in the space $C(\A)$.\\
\vspace{-0.2cm}
\\
$\small{\textbf{\quad b) Finite dimensional convergence}}$\\
\vspace{-0.2cm}
\\
The convergence of the finite-dimensional laws is a simple
consequence of both the central limit theorem for random fields
(\cite{JD-tcl}, Theorem 2.2) and the following lemma (see
\cite{JD-tflc}). For any subset $\Gamma$ of $\Z^{d}$ we consider
$$
\partial\Gamma=\big\{i\in\Gamma\,;\,\exists{j\notin\Gamma}\,
\,\textrm{such that}\,\,\vert i-j\vert=1\big\}.
$$
For any Borel set $A$ of $[0,1]^{d}$, we denote by $\Gamma_{n}(A)$
the finite subset of $\Z^{d}$ defined by
$\Gamma_{n}(A)=nA\cap\Z^{d}$.
\begin{Lem}[Dedecker, 2001]\label{lemma2}
Let $A$ be a regular Borel set of $[0,1]^{d}$ with $\lambda(A)>0$.
We have
$$
(i)\,\,\lim_{n\to+\infty}\frac{\vert\Gamma_{n}(A)\vert}{n^{d}}=\lambda(A)\qquad
(ii)\,\,\lim_{n\to+\infty}\frac{\vert\partial\Gamma_{n}(A)\vert}{\vert\Gamma_{n}(A)\vert}=0.
$$
Let $(X_{i})_{i\in\Z^{d}}$ be a stationary random field with mean
zero and finite variance.\\
Assume that $\sum_{k\in\Z^{d}}\vert E(X_{0}X_{k})\vert<+\infty$.
Then
\begin{displaymath}
\lim_{n\to+\infty}\quad
n^{-d/2}\bigg\|S_{n}(A)-\sum_{k\in\Gamma_{n}(A)}X_{k}\bigg\|_{2}=0.
\end{displaymath}
\end{Lem}
\subsection{Proof of Theorem $\bold{\ref{fclt}}$}
Similarly, we are going to prove both the convergence of the
finite-dimensional laws and the tightness of the sequence of
processes $\{U_{n}^{-1}S_{n}(A)\,;\,A\in\A\}$ in the space
$C(\A)$.\\
\vspace{-0.2cm}
\\
$\small{\textbf{\quad a) Tightness}}$\\
\vspace{-0.2cm}
\\
It suffices to establish that for any $x>0$
\begin{equation}\label{tension}
\lim_{\delta\to0}\limsup_{n\to+\infty}
\P\left(\sup_{\substack{A,B\in\A \\\rho(A,B)<\delta}}\big\vert
U_{n}^{-1}S_{n}(A)-U_{n}^{-1}S_{n}(B)\big\vert>x\right)=0.
\end{equation}
Let $\delta>0$ and $0<\tau\leq 1$ defined as in the proof of
theorem $\ref{fclt-am}$. In the sequel, we denote $(b_{n})_{n\geq
1}$ the sequence which satisfies condition ($\ref{standard}$) and
we define $X_{i,n}=X_{i}\ind{{\vert X_{i}\vert<\tau b_{n}}}$. One
can check that
$$
\P\left(\sup_{\substack{A,B\in\A \\ \rho(A,B)<\delta}}\big\vert
U_{n}^{-1}S_{n}(A)-U_{n}^{-1}S_{n}(B)\big\vert>x\right)\leq
E_{1}+E_{2}+E_{3}+E_{4}
$$
where
\begin{align*}
E_{1}&=\P\left(\sup_{\substack{A,B\in\A \\ \rho(A,B)<\delta}}\bigg\vert\sum_{i\in\Lambda{n}}(\lambda(nA\cap R_{i})-\lambda(nB\cap R_{i}))[X_{i,n}-EX_{i,n}]\bigg\vert>xb_{n}/2\right)\\
E_{2}&=\P\left(U_{n}\leq b_{n}/2\right)\tend{n}{+\infty}0\qquad\textrm{(by Raikov's theorem)}\\
E_{3}&=n^d\P\left(\vert X_{0}\vert\geq\tau b_{n}\right)\tend{n}{+\infty}0\qquad\textrm{(by ($\ref{standard}$))}\\
E_{4}&=x^{-1}b_{n}^{-1}n^d \vert
EX_{0,n}\vert\tend{n}{+\infty}0\qquad\textrm{(by
($\ref{standard}$))}.
\end{align*}
So, it suffices to control $E_{1}$. As in the proof of Theorem
$\ref{fclt-am}$, we apply the chaining method by Bass \cite{Bass}
with the following notations: for any constants
$0\leq\alpha\leq\beta\leq 1$, we define
$X_{i}(n,\alpha,\beta)=X_{i}\ind{\alpha\tau b_{n}\leq\vert
X_{0}\vert<\beta\tau b_{n}}$ and
$$
Z_{n}(A,\alpha,\beta)=\frac{1}{b_{n}}\sum_{i\in\Lambda_{n}}\lambda(nA\cap
R_{i})[X_{i}(n,\alpha,\beta)-EX_{i}(n,\alpha,\beta)].
$$
So, we obtain
$$
E_{1}\leq\frac{2}{x}E\left(\sup_{\substack{A,B\in\A \\
\rho(A,B)<\delta}}\big\vert
Z_{n}(A,0,1)-Z_{n}(B,0,1)\big\vert\right)\leq\frac{2}{x}\left(F_{1}+F_{2}+F_{3}\right)
$$
where $F_{1}$, $F_{2}$ and $F_{3}$ are defined in the proof of
Theorem $\ref{fclt-am}$. Applying Lemma $\ref{Lem1}$ with
$c_{n}=b_{n}$, the estimations ($\ref{F1}$) and ($\ref{F2}$) still
hold for $F_{1}$ and $F_{2}$ respectively. In order to control the
last term $F_{3}$, for any Borel subset $A$ of $[0,1]^d$, we
denote
$$
\widetilde{Z}_{n}(A,a_{k},a_{k-1})=\frac{1}{b_{n}}\sum_{i\in\Lambda_{n}}\lambda(nA\cap
R_{i})[\vert X_{i}(n,a_{k},a_{k-1})\vert-E\vert
X_{i}(n,a_{k},a_{k-1})\vert].
$$
We have
\begin{align*}
&\sup_{A_{k}\subset A\subset A_{k}^{+}}\vert Z_{n}(A,a_{k},a_{k-1})-Z_{n}(A_{k},a_{k},a_{k-1})\vert\\
&\qquad\leq
\frac{1}{b_{n}}\sum_{i\in\Lambda_{n}}(\lambda(nA_{k}^{+}\cap
R_{i})-\lambda(nA_{k}\cap R_{i}))
[\vert X_{i}(n,a_{k},a_{k-1})\vert-E\vert X_{i}(n,a_{k},a_{k-1})\vert]\\
&\qquad\qquad+2\frac{n^d}{b_{n}}E\vert
X_{0}(n,a_{k},a_{k-1})\vert\,\delta_{k}^2\\
&\qquad=\widetilde{Z}_{n}(A_{k}^{+},a_{k},a_{k-1})-\widetilde{Z}_{n}(A_{k},a_{k},a_{k-1})+2\frac{n^d}{b_{n}}E\vert
X_{0}(n,a_{k},a_{k-1})\vert\,\delta_{k}^2
\end{align*}
Using Lemma $\ref{Lem1}$, we derive
$$
\bigg\|\max_{A_{k},A_{k}^{+}\in\A(\delta_{k})}\big\vert\widetilde{Z}_{n}(A_{k}^{+},a_{k},a_{k-1})-\widetilde{Z}_{n}(A_{k},a_{k},a_{k-1})\big\vert\bigg\|_{\psi_{1}}
\hspace{-0.2cm}\leq
K(a_{k-1}\tau\H(\delta_{k})+\delta_{k}\sqrt{\H(\delta_{k})}).
$$
In the other hand
$$
\frac{n^d}{b_{n}}E\vert X_{0}(n,a_{k},a_{k-1})\vert\,\delta_{k}^2
\leq\frac{\delta_{k}^2}{a_{k}\tau}\frac{n^d}{b_{n}^2}E
X_{0}^2\ind{{\vert X_{0}\vert<b_{n}}}.
$$
So, the estimation ($\ref{F3}$) still hold for $F_{3}$ and
choosing again $a_{k}=\delta_{k}/(\tau\sqrt{\H(\delta_{k+1})})$,
we derive
$$ \limsup_{n\to\infty}\frac{2}{x}E_{1}\leq
K\sum_{k=1}^{+\infty}\delta_{k+1}\sqrt{\H(\delta_{k})}\leq
K\int_{0}^{\delta}\sqrt{\H(x)}dx\tend{\delta}{0}0.
$$
Finally, the condition ($\ref{tension}$) holds and the sequence
$\{U_{n}^{-1}S_{n}(A)\,;\,A\in\A\}$ is tight in the
space $C(\A)$.\\
\vspace{-0.2cm}
\\
$\small{\textbf{\quad b) Finite dimensional convergence}}$\\
\vspace{-0.2cm}
\\
For any Borel set $A$ of $[0,1]^{d}$ recall that $\Gamma_{n}(A)$
is the finite set defined by $\Gamma_{n}(A)=nA\cap\Z^{d}$ and
denote $S_{\Gamma_{n}(A)}=\sum_{i\in\Gamma_{n}(A)}X_{i}$.
\begin{Lem}\label{finidi1}
Let $A$ be a regular Borel set of $[0,1]^{d}$ with $\lambda(A)>0$.
For any $x>0$, we have
$$
\lim_{n\to\infty}\P\left( U_{n}^{-1}\vert
S_{n}(A)-S_{\Gamma_{n}(A)}\vert>x\right)=0.
$$
\end{Lem}
{\em Proof of Lemma $\ref{finidi1}$}. Consider the subsets of
$\Z^{d}$
$$
A_{1}=\{i\,;\,R_{i}\subset nA\},\quad A_{2}=\{i\,;\,R_{i}\cap
nA\neq\emptyset\},\quad A_{3}=A_{2}\cap \{i\,;\,R_{i}\cap
(nA)^{c}\neq\emptyset\}
$$
and set $a_{i}=\lambda(nA\cap R_{i})-\ind{i\in\Gamma_{n}(A)}$.
Since $a_{i}$ equals zero if $i$ belongs to $A_{1}$, we have
$$
S_{n}(A)-S_{\Gamma_{n}(A)}=\sum_{i\in A_{3}}a_{i}X_{i}.
$$
Let $\tau>0$ and recall that $X_{i,n}=X_{i}\ind{\vert
X_{i}\vert<\tau b_{n}}$. We have
$$
\P\left( U_{n}^{-1}\vert
S_{n}(A)-S_{\Gamma_{n}(A)}\vert>x\right)\leq P_{1}+P_{2}+P_{3}
$$
where
\begin{align*}
P_{1}&=\P\left(\bigg\vert\sum_{i\in A_{3}}a_{i}X_{i,n}\bigg\vert>xb_{n}/2\right)\\
P_{2}&=\P\left(U_{n}\leq
b_{n}/2\right)\tend{n}{+\infty}0\quad\textrm{(by ($\ref{Raikov}$))}\\
P_{3}&=n^d\P\left(\vert X_{0}\vert\geq\tau
b_{n}\right)\tend{n}{+\infty}0\quad\textrm{(by
($\ref{standard}$))}.
\end{align*}
Moreover
$$
P_{1}\leq\frac{4\vert A_{3}\vert}{x^2b_{n}^2}
EX_{0,n}^2=\frac{4\vert
A_{3}\vert}{x^2n^d}\times\frac{n^d}{b_{n}^2}EX_{0,n}^2.
$$
Keeping in mind that $n^{-d}\vert A_{3}\vert$ tends to zero as $n$
goes to infinity (cf. Dedecker \cite{JD-tflc}) and using
($\ref{standard}$) then the proof of Lemma $\ref{finidi1}$ is
complete.
\begin{Lem}\label{finidi2}
For any regular Borel set $A$ in $\A$, the sequence
$\left(U_{n}^{-1}S_{\Gamma_{n}(A)}\right)_{n\geq 1}$ converge in
distribution to $\sqrt{\lambda(A)}\,\varepsilon$ where
$\varepsilon$ has the standard normal law.
\end{Lem}
{\em Proof of Lemma $\ref{finidi2}$}. Let $x>0$, $n\in\N^{\ast}$
and $A\in\A$ be fixed. We have
$$
U_{n}^{-1}S_{\Gamma_{n}(A)}=\underbrace{\frac{\sum_{i\in\Gamma_{n}(A)}X_{i}}{\sqrt{\sum_{i\in\Gamma_{n}(A)}X_{i}^{2}}}}_{T_{n,1}(A)}
\times\underbrace{\sqrt{\frac{\sum_{i\in\Gamma_{n}(A)}X_{i}^2}{\sum_{i\in\Lambda_{n}}X_{i}^2}}}_{T_{n,2}(A)}.
$$
Using Theorem 3.3 in \cite{Gine-Gotze-Mason}, we derive that
$T_{n,1}(A)$ converges in distribution to the standard normal law.
So, it suffices to prove that $T_{n,2}^2(A)$ converges in
probability to $\lambda(A)$. Let $\tau>0$ be fixed. Denoting
$X_{i,n}=X_{i}\ind{\vert X_{i}\vert<\tau b_{n}}$ and
$\overline{X}_{i,n}=X_{i}-X_{i,n}$, we have
\begin{equation}\label{ineq-triang}
\vert T_{n,2}^2(A)-\lambda(A)\vert \leq\underbrace{\bigg\vert
T_{n,2}^2(A)-\frac{\sum_{i\in\Gamma_{n}(A)}X_{i,n}^2}{\sum_{i\in\Lambda_{n}}X_{i,n}^2}\bigg\vert}_{(\ast)}
+\underbrace{\bigg\vert\frac{\sum_{i\in\Gamma_{n}(A)}X_{i,n}^2}{\sum_{i\in\Lambda_{n}}X_{i,n}^2}-\lambda(A)\bigg\vert}_{(\ast\ast)}.
\end{equation}
Now, noting that $X_{i}^2=X_{i,n}^2+\overline{X}_{i,n}^2$, we derive
\begin{align*}
(\ast)&=\bigg\vert\frac{\sum_{i\in\Lambda_{n}}X_{i,n}^2\sum_{i\in\Gamma_{n}(A)}X_{i}^2
-\sum_{i\in\Lambda_{n}}X_{i}^2\sum_{i\in\Gamma_{n}(A)}X_{i,n}^2}{\sum_{i\in\Lambda_{n}}X_{i}^2\sum_{i\in\Lambda_{n}}X_{i,n}^2}\bigg\vert\\
&=\bigg\vert\frac{\sum_{i\in\Lambda_{n}}X_{i,n}^2\sum_{i\in\Gamma_{n}(A)}\overline{X}_{i,n}^2
-\sum_{i\in\Lambda_{n}}\overline{X}_{i,n}^2\sum_{i\in\Gamma_{n}(A)}X_{i,n}^2}{\sum_{i\in\Lambda_{n}}X_{i}^2\sum_{i\in\Lambda_{n}}X_{i,n}^2}\bigg\vert\\
&\leq
2\,\frac{\sum_{i\in\Lambda_{n}}\overline{X}_{i,n}^2}{\sum_{i\in\Lambda_{n}}X_{i}^2}\\
&=2\left(1-R_{n}\right)
\end{align*}
where
$$
R_{n}=\frac{\sum_{i\in\Lambda_{n}}X_{i,n}^2}{\sum_{i\in\Lambda_{n}}X_{i}^2}\leq
1\quad\textrm{a.s.}
$$
Let $x>0$ be fixed. Using ($\ref{standard}$) we derive that
\begin{equation}\label{ast}
\P((\ast)>3x)\leq\P((\ast)>0)\leq\P(R_{n}<1)\leq
n^d\P(\vert X_{0}\vert\geq\tau b_{n})\tend{n}{+\infty}0.
\end{equation}
In the other hand,
\begin{align*}
(\ast\ast)&\leq\bigg\vert\frac{\sum_{i\in\Gamma_{n}(A)}X_{i,n}^2}{\sum_{i\in\Lambda_{n}}X_{i,n}^2}-
\frac{1}{b_{n}^2}\sum_{i\in\Gamma_{n}(A)}X_{i,n}^2\bigg\vert+
\bigg\vert\frac{1}{b_{n}^2}\sum_{i\in\Gamma_{n}(A)}X_{i,n}^2-\lambda(A)\bigg\vert\\
&\leq\bigg\vert
1-\frac{1}{b_{n}^2}\sum_{i\in\Lambda_{n}}X_{i,n}^2\bigg\vert
+\bigg\vert\frac{1}{b_{n}^2}\sum_{i\in\Gamma_{n}(A)}X_{i,n}^2-\lambda(A)\bigg\vert\\
&\leq\underbrace{\bigg\vert
1-\frac{1}{b_{n}^2}\sum_{i\in\Lambda_{n}}X_{i,n}^2\bigg\vert}_{\gamma_{n,1}}
+\underbrace{\bigg\vert\frac{1}{b_{n}^2}\sum_{i\in\Gamma_{n}(A)}\left(X_{i,n}^2-EX_{i,n}^2\right)\bigg\vert}_{\gamma_{n,2}}
+\underbrace{\bigg\vert\frac{\vert\Gamma_{n}(A)\vert}{b_{n}^2}EX_{0,n}^2-\lambda(A)\bigg\vert}_{\gamma_{n,3}}.
\end{align*}
By ($\ref{standard}$) and the point $(i)$ of Lemma $\ref{lemma2}$,
it is clear that
\begin{equation}\label{gamma-n3}
\gamma_{n,3}\tend{n}{\infty}0.
\end{equation}
Noting that
$$
b_{n}^{-2}\sum_{i\in\Lambda_{n}}X_{i,n}^2=\frac{\sum_{i\in\Lambda_{n}}X_{i}^2}{b_{n}^2}\times
R_{n}\quad\textrm{a.s.}
$$
we have
\begin{align*}
\P(\gamma_{n,1}>x)&\leq\P\left(\big\vert
1-R_{n}\big\vert>x/2\right)+\P\left(\bigg\vert
1-\frac{\sum_{i\in\Lambda_{n}}X_{i}^2}{b_{n}^2}\bigg\vert>x/2\right)\\
&\leq\P(R_{n}<1)+\P\left(\bigg\vert
1-\frac{\sum_{i\in\Lambda_{n}}X_{i}^2}{b_{n}^2}\bigg\vert>x/2\right)\\
&\leq n^d\P(\vert X_{0}\vert\geq\tau b_{n})+\P\left(\bigg\vert
1-\frac{\sum_{i\in\Lambda_{n}}X_{i}^2}{b_{n}^2}\bigg\vert>x/2\right).
\end{align*}
Using ($\ref{standard}$) and ($\ref{Raikov}$), we obtain
\begin{equation}\label{gamma-n1}
\P(\gamma_{n,1}>x)\tend{n}{\infty}0.
\end{equation}
We have also
\begin{align*}
\P(\gamma_{n,2}>x)&\leq\frac{b_{n}^{-4}}{x^2}E\left(\sum_{i\in\Gamma_{n}(A)}X_{i,n}^2-EX_{i,n}^2\right)^2\\
&=\frac{b_{n}^{-4}}{x^2}\vert\Gamma_{n}(A)\vert
E\left(X_{0,n}^2-EX_{0,n}^2\right)^2\\
&\leq\frac{4b_{n}^{-4}}{x^2}\vert\Gamma_{n}(A)\vert EX_{0,n}^4\\
&\leq\frac{4\tau^2b_{n}^{-2}}{x^2}\vert\Gamma_{n}(A)\vert EX_{0,n}^2\\
&=\frac{4\tau^2\vert\Gamma_{n}(A)\vert}{n^d
x^2}\times\frac{n^d}{b_{n}^2}EX_{0,n}^2.
\end{align*}
Consequently, using ($\ref{standard}$) and the point $(i)$ in
Lemma $\ref{lemma2}$, we derive
\begin{equation}\label{gamma-n2}
\lim_{n\to+\infty}\P(\gamma_{n,2}>x)\leq\frac{4\tau^2\lambda(A)}{x^2}.
\end{equation}
Now, combining ($\ref{gamma-n3}$), ($\ref{gamma-n1}$) and ($\ref{gamma-n2}$), we obtain
\begin{equation}\label{ast-ast}
\lim_{n\to+\infty}\P((\ast\ast)>3x)\leq\frac{4\tau^2\lambda(A)}{x^2}.
\end{equation}
Combining ($\ref{ineq-triang}$), ($\ref{ast}$)
and ($\ref{ast-ast}$), it follows that
$$
\lim_{n\to+\infty}\P\left(\vert T_{n,2}^2(A)-\lambda(A)\vert>6x\right)\leq\frac{4\tau^2\lambda(A)}{x^2}.
$$
Since $\tau>0$ can be arbitrarily small, we obtain
$$
\lim_{n\to+\infty}\P\left(\vert T_{n,2}^2(A)-\lambda(A)\vert>6x\right)=0.
$$
Finally, $T_{n,2}^2(A)$ converges in probability to $\lambda(A)$
and the proof of Lemma $\ref{finidi2}$ is complete. The
convergence of the finite-dimensional laws of the sequence
$\{U_{n}^{-1}S_{n}(A);A\in\A\}$ follows then from Lemmas
$\ref{finidi1}$ and $\ref{finidi2}$. The proof of Theorem
$\ref{fclt}$ is complete.
\subsection{Proof of Theorem $\bold{\ref{counter}}$}
Without loss of generality, we assume that $p$ is a positive
integer. Consider the field $X=(X_{k})_{k\in\Z^{d}}$ of i.i.d.
integer-valued random variables defined on a probability space
$(\Omega, \F, \mu)$ by the following property: the random variable
$X_{0}$ is symmetric and satisfies $\mu(X_{0}=0)=0$ and $\mu(\vert
X_{0}\vert\geq k)=k^{-p-1}$ for any integer $k\geq 1$. The random
field $X$ is $p$-integrable since
\begin{align*}
E(\vert X_{0}\vert^{p})&=\sum_{k\geq 1}\mu(\vert X_{0}\vert\geq k^{1/p})\\
&=\sum_{k\geq 1}k^{-1-1/p}<+\infty.
\end{align*}
Let us fix an integer $r\geq 1$ and consider the following
numbers:
$$
n_{r}=4^{rp},
$$
$$
\beta_{r}=n_{r}^{d/2p}=2^{rd},
$$
$$
k_{r}=n_{r}^{d}\mu(X_{0}\geq\beta_{r})=2^{rd(p-1)},
$$
$$
\varepsilon_{r}=\left(\frac{k_{r}}{n_{r}^{d}}\right)^{1/2}=2^{-rd(p+1)/2}.
$$
One can notice that $(n_{r})_{r\geq 1}$, $(\beta_{r})_{r\geq 1}$
and $(k_{r})_{r\geq 1}$ are increasing sequences of positive
integers while $(\varepsilon_{r})_{r\geq 1}$ is a decreasing
sequence of positive real numbers which converges to zero. We
define the class $\A_{r}$ as the collection of all Borel subsets
$A$ of $[0,1]^{d}$ with the following property: $A$ is empty or
there exist $i_{l}=(i_{l,1},...,i_{l,d})$ in
$\{1,...,n_{r}\}^{d},\,1\leq l\leq k_{r}$ such that
$$
A=\ds{\bigcup_{l=1}^{k_{r}}}\,\,\,
\bigg]\frac{i_{l,1}-1}{n_{r}},\frac{i_{l,1}}{n_{r}}\bigg]\times...\times
\bigg]\frac{i_{l,d}-1}{n_{r}},\frac{i_{l,d}}{n_{r}}\bigg].
$$
Now, denote
$$
\A=\B_{r}\cup\C_{r}
$$
where
$$
\B_{r}=\bigcup_{j=1}^{r-1}\A_{j}\quad\textrm{and}\quad\C_{r}=\bigcup_{j=r}^{+\infty}\A_{j}.
$$
For any integer $j\geq 1$, the cardinal $\vert\A_{j}\vert$ of $\A_{j}$ equals
$1+\left(\begin{array}[c]{c}n_{j}^{d}\\k_{j} \end{array}\right)$,
hence
$$
N(\B_{r}, \rho,
\varepsilon_{r})\leq\sum_{j=1}^{r-1}\left(1+\left(\begin{array}[c]{c}n_{j}^{d}\\k_{j}
\end{array}\right)\right)\leq 2rn_{r}^{dk_{r}}.
$$
On the other hand, since each element of the class $\C_{r}$
belongs to the ball with center $\emptyset$ and radius
$\varepsilon_{r}$, it follows that $N(\C_{r}, \rho,
\varepsilon_{r})=1$. Noting that
$$
N(\A, \rho, \varepsilon_{r}) \leq N(\B_{r}, \rho,
\varepsilon_{r})+N(\C_{r}, \rho, \varepsilon_{r}),
$$
we obtain
$$
N(\A, \rho, \varepsilon_{r})\leq 1+2rn_{r}^{dk_{r}}
$$
and also
$$
H(\A, \rho, \varepsilon_{r})=\log\,N(\A, \rho,
\varepsilon_{r})\leq 3dk_{r}\log\,n_{r}.
$$
Finally, there exists $K>0$ such that
\begin{align*}
\sum_{r=2}^{+\infty}\varepsilon_{r-1}\sqrt{H(\A, \rho,
\varepsilon_{r})}
&\leq \sum_{r=2}^{+\infty}\varepsilon_{r-1}\sqrt{3dk_{r}\log\,n_{r}}\\
&\leq K
\sum_{r=2}^{+\infty}\frac{2^{rd(p-1)/2}\sqrt{r}}{2^{rd(p+1)/2}}\\
&=K\sum_{r=2}^{+\infty}\frac{\sqrt{r}}{2^{rd}}<+\infty.
\end{align*}
Consequently, the class $\A$ satisfies the metric entropy
condition ($\ref{entrop-metriq1}$). Now, we are going to see that
the partial sum process $\{n^{-d/2}S_{n}(A)\,;\,A\in\A\}$ defined
by ($\ref{process}$) is not tight in the space $C(\A)$. It is
sufficient (Pollard, 1990\nocite{Pollard}) to show that there
exists $\theta>0$ such that
$$
\lim_{\delta\to 0}\limsup_{n\to+\infty}\mu\left(\sup_{\substack{A,B\in\A \\
\rho(A,B)<\delta}} n^{-d/2}\big\vert
S_{n}(A)-S_{n}(B)\big\vert\geq\theta\right)>0.
$$
For any integer $r\geq 1$, denote
$\Lambda_{r}=\{1,...,n_{r}\}^{d}$ and define $W_{r}$ as the set of
all $\omega$ in $\Omega$ such that
$$
\sum_{i\in\Lambda_{r}}\ind{\{X_{i}(\omega)\geq\beta_{r}\}}\geq
k_{r}.
$$
\begin{Lem}\label{lemme-minoration}
There exists a constant $c>0$ such that for any integer $r\geq 1$,
\begin{equation}\label{minoration}
\mu(W_{r})\geq c.
\end{equation}
\end{Lem}
{\em Proof of Lemma $\ref{lemme-minoration}$}. Let $r\geq 1$ be fixed. For any $i$
in $\Lambda_{r}$, denote
$$
Y_{i}=\ind{\{X_{i}\geq\beta_{r}\}}-\mu(X_{0}\geq\beta_{r}).
$$
The family $\{Y_{i}\,;\,i\in\Lambda_{r}\}$ is a finite sequence of
i.i.d. centered random variables bounded by 2. So, using a lower
exponential inequality due to Kolmogorov (Ledoux and Talagrand, 1991\nocite{Led-Tal},
Lemma 8.1), it follows that for any $\gamma>0$, there exist
positive numbers $K(\gamma)$ (large enough) and
$\varepsilon(\gamma)$ (small enough) depending on $\gamma$ only,
such that for every $t$ satisfying $t\geq K(\gamma)b$ and
$2t\leq\varepsilon(\gamma)b^{2}$,
$$
\mu\left(\sum_{i\in\Lambda_{r}}Y_{i}>t\right)\geq\exp\left(-(1+\gamma)t^{2}/2b^{2}\right)
$$
where $b^{2}=\sum_{i\in\Lambda_{r}}EY_{i}^{2}$. In particular,
there exists a positive universal constant $K$ such that
$$
\mu\left(\sum_{i\in\Lambda_{r}}Y_{i}>Kb\right)\geq\exp\left(-K^{2}\right).
$$
Noting $c=\exp(-K^{2})>0$ and keeping in mind the definitions of
the constant $k_{r}$ and the random variable $Y_{i}$, we derive
$$
\mu\left(\sum_{i\in\Lambda_{r}}\ind{\{X_{i}\geq\beta_{r}\}}>Kb+k_{r}
\right)\geq c.
$$
Finally, Inequality ($\ref{minoration}$) follows from the fact
that $Kb\geq 0$ and the proof of the lemma is complete. The proof of Lemma $\ref{lemme-minoration}$ is complete.\\
\\
Let $\omega$ be fixed in the set $W_{r}$ and denote
$$
\Gamma_{r}^{\ast}(\omega)=\{i\in\Lambda_{r}\,;\,X_{i}(\omega)\geq\beta_{r}\}.
$$
By definition of the set $W_{r}$, we know that
$\vert\Gamma_{r}^{\ast}(\omega)\vert\geq k_{r}$. Let
$\Gamma_{r}(\omega)$ be a subset of $\Gamma_{r}^{\ast}(\omega)$
such that $\vert\Gamma_{r}(\omega)\vert=k_{r}$ and define
$$
A_{r}(\omega)=\ds{\bigcup_{i\in\Gamma_{r}(\omega)}}\,\,\,
\bigg]\frac{i_{1}-1}{n_{r}},\frac{i_{1}}{n_{r}}\bigg]\times...\times
\bigg]\frac{i_{d}-1}{n_{r}},\frac{i_{d}}{n_{r}}\bigg]\in\A_{r}\subset\A.
$$
For any $\omega$ in $W_{r}$ and any $i$ in $\Lambda_{r}$, we have
$$
\lambda(n_{r}A_{r}(\omega)\cap R_{i})=\ind{\Gamma_{r}(\omega)}(i).
$$
Consequently, we have
\begin{align*}
n_{r}^{-d/2}S_{n_{r}}(A_{r}(\omega))&=n_{r}^{-d/2}\sum_{i\in\Lambda_{r}}\lambda(n_{r}A_{r}(\omega)\cap
R_{i})X_{i}(\omega)\\
&=n_{r}^{-d/2}\sum_{i\in\Gamma_{r}(\omega)}X_{i}(\omega)\\
&\geq n_{r}^{-d/2}\vert\Gamma_{r}(\omega)\vert\beta_{r}\\
&=n_{r}^{-d/2}k_{r}\beta_{r}\\
&=n_{r}^{d/2}\mu(X_{0}\geq\beta_{r})\beta_{r}\\
&=\frac{1}{2}n_{r}^{d/2}\beta_{r}^{-p}\\
&=\frac{1}{2}.
\end{align*}
Thus, for any integer $r\geq 1$ and any $\omega$ in $W_{r}$, we
have
\begin{equation}\label{minoration2}
\big\vert n_{r}^{-d/2}S_{n_{r}}(A_{r}(\omega))\big\vert\geq 1/2.
\end{equation}
Let $\delta>0$ be fixed. There exists an integer $R$ such that for
any $r\geq R$ and any $\omega$ in $W_{r}$,
$\lambda(A_{r}(\omega))=k_{r}/n_{r}^{d}\leq\delta^{2}$. Then,
using the lower bounds ($\ref{minoration}$) and
($\ref{minoration2}$), it follows that for any $r\geq R$,
\begin{align*}
&\mu\left(\sup_{\substack{A,B\in\A
\\ \rho(A,B)<\delta}}\big\vert
n_{r}^{-d/2}S_{n_{r}}(A)-n_{r}^{-d/2}S_{n_{r}}(B)\big\vert\geq
1/2\right)\\
&\geq \mu\left(\sup_{\substack{A\in\A
\\ \lambda(A)<\delta^{2}}}\big\vert
n_{r}^{-d/2}S_{n_{r}}(A)\big\vert\geq
1/2\right)\\
&\geq \mu\left(\bigg\{\omega\in W_{r}\,\,\bigg\vert\,\,\big\vert
n_{r}^{-d/2}S_{n_{r}}(A_{r}(\omega))\big\vert\geq
1/2\bigg\}\right)\\
&=\mu(W_{r})\geq c>0.
\end{align*}
Finally, we have shown that for any $\delta>0$,
$$
\limsup_{n\to+\infty}\mu\left(\sup_{\substack{A,B\in\A
\\ \rho(A,B)<\delta}}\big\vert
n^{-d/2}S_{n}(A)-n^{-d/2}S_{n}(B)\big\vert\geq 1/2\right)\geq c>0.
$$
The proof of Theorem $\ref{counter}$ is complete.\\
\\
\textbf{Acknowledgement}. The authors thank Prof. Dalibor Voln\'y
for his useful help in the construction of the counter-example in
Theorem $\ref{counter}$.
\bibliographystyle{plain}
\bibliography{xbib}
\end{document}